\documentclass{amsart}

\usepackage[ngerman,american]{babel}
\usepackage{a4wide}
\usepackage{amssymb}
\usepackage[latin1]{inputenc}

\renewcommand{\phi}{\varphi}
\newcommand{\C}{{\mathbb{C}}}
\newcommand{\R}{{\mathbb{R}}}

\newcommand{\N}{{\mathbb{N}}}

\renewcommand{\epsilon}{\varepsilon}
\renewcommand{\theta}{\vartheta}
\newcommand{\D}{D}
\renewcommand{\S}{S}

\newcommand{\Brieskorn}[2]{{W^{#1}_{#2}}}

\newcommand{\pairing}[2]{{\langle{#1}|{#2}\rangle}}
\newcommand{\lcan}{{\lambda_{\mathrm{can}}}}

\newcommand{\Dehn}{\tau}
\newcommand{\POS}{{\mathbf{q}}}
\newcommand{\MOM}{{\mathbf{p}}}

\DeclareMathOperator{\id}{id}
\DeclareMathOperator{\orb}{Orb}

\DeclareMathOperator{\SO}{SO}

\theoremstyle{plain}

\theoremstyle{remark}
\newtheorem*{remark}{Remark}

\theoremstyle{definition}
\newtheorem{defi}{Definition}

\begin{document}
\bibliographystyle{amsalpha}

\title[Open books]{Open book decompositions for contact structures on
  Brieskorn manifolds}
\author{Otto van Koert}
\email{okoert@math.uni-koeln.de}
\author{Klaus Niederkrüger}
\email{kniederk@math.uni-koeln.de}
\address{Mathematisches Institut, Universität zu Köln\\Weyertal 86-90\\50.931 Köln\\Federal Republic of Germany}

\begin{abstract}
  In this paper, we give an open book decomposition for the contact
  structures on some Brieskorn manifolds, in particular for the contact
  structures of Ustilovsky.
  The decomposition uses right-handed Dehn twists as conjectured by Giroux.
\end{abstract}

\maketitle

\setcounter{section}{-1}
\section{Introduction}
At the ICM of 2002 Giroux announced some of his results concerning a
correspondence between contact structures on manifolds and open book
structures on them. In one direction this correspondence is relatively easy. We are
given a compact Stein manifold $P$ (i.e.~a compact subset of a Stein manifold
where the boundary is a level set of a plurisubharmonic function on it) and
a symplectomorphism $\psi$ of $P$ that is the identity near the boundary of
$P$. It can be shown that this symplectomorphism gives rise to a mapping
torus that inherits a contact structure. Furthermore
the boundary of the mapping torus will always look like $\S^1 \times
\partial P$, and the binding, $\D^2 \times \partial P$ with the obvious contact
structure, can be glued in to give a compact contact manifold.

Although Giroux announced much more than just this, it is already
interesting to see how this construction turns out in a few simple cases. As a
Stein manifold we will take $T^*\S^{n-1}$ with its canonical symplectic
form. The symplectomorphisms used for the monodromy of the mapping
torus will be so-called generalized Dehn twists. Seidel has shown \cite{Seidel_Preprint}
that these Dehn twists generate the symplectomorphism group of $T^*\S^2$ up
to isotopy. Furthermore his results show that Dehn twists of $T^*\S^2$ are
of order 2 diffeomorphically, but not symplectically. This means that many
of these Dehn twists are isotopic to each other, but not symplectically so.  

In the spirit of the above construction, we will show that the Brieskorn
manifold $W_k^{2n-1}$ (for notation and definition see Section~\ref{brieskorn}) is supported by an open book whose monodromy is given by
a $k$-fold Dehn twist. In particular this shows that the Ustilovsky
spheres (special Brieskorn spheres with non-isomorphic contact structures)
can all be written in terms of open book decompositions with Dehn twists as their
monodromy. It also shows that Dehn twists cannot be of order 2 in all
dimensions (this is well known for $n$ even). Namely, among the Brieskorn
spheres (these correspond to $n$ and $k$ odd) are exotic spheres as well as
standard ones. As the binding is always glued in in the same way, the Dehn
twists corresponding to a standard and an exotic sphere cannot be isotopic
relative to the boundary.

\subsubsection*{Acknowledgements}
The authors would like to thank F.\ Bourgeois and H.\ Geiges for
their helpful comments on this paper.

\section{Notation \& Definitions}\label{notation}

\subsection{Open books}

The following definitions are taken from \cite{Giroux}.

\begin{defi}
  An \textbf{open book} on a closed manifold $M$ is given by a codimension-$2$ submanifold $B\hookrightarrow M$ with trivial normal bundle, and a bundle $ \theta:\,(M-B)\to \S^1$.
  The neighborhood of $B$ should have a trivialization $B\times\D^2$, where the angle coordinate on the disk agrees with the map $\theta$.

  The manifold $B$ is called the \textbf{binding} of the open book and a fiber $P=\theta^ {-1}(\varphi_0)$ is called a \textbf{page}.
\end{defi}

\begin{remark}
  The open set $M-B$ is a bundle over $\S^1$, hence it is diffeomorphic to $\R\times P/\sim$, where $\sim$ identifies $(t,p)\sim(t+1,\Phi(p))$ for some diffeomorphism $\Phi$ of $P$.
\end{remark}

\begin{defi}\label{compatiblebook}
A contact structure $\xi=\ker\alpha$ on $M$ is said to be \textbf{supported
  by an open book $(B,\theta)$} of $M$, if
\begin{enumerate}
\item $(B,\alpha|_{TB})$ is a contact manifold.
\item For every $s\in \S^1$, the page $P:=\theta^{-1}(s)$ is a
  symplectic manifold with symplectic form $d\alpha$.
\item Denote the closure of a page $P$ in $M$ by $\bar P$. The orientation
  of $B$ induced by its contact form $\alpha|_{TB}$ should coincide
  with its orientation as the boundary of $(\bar P,d\alpha)$.
 
\end{enumerate} 
Such a contact form is said to be \textbf{adapted} to $(B,\theta)$.
\end{defi}

\subsection{Dehn twists}\label{dehn twist}
A \textbf{Dehn twist $\Dehn_k$} is a diffeomorphism from $T^*\S^{n-1}$ to itself constructed in the following way.
Write points in $T^*\S^{n-1}$ as $(\POS,\MOM)\in\R^{2n}$ with $|\POS|=1$ and $\POS\perp\MOM$.

Set
$$
\tau_k(\POS,\MOM)=
\left(
\begin{array}{cc}
\cos g_k(\MOM) & |\MOM|^{-1} \sin g_k(\MOM) \\
-|\MOM|\sin g_k(\MOM) & \cos g_k(\MOM) 
\end{array}
\right)
\left(
\begin{array}{c}
\POS \\ \MOM
\end{array}
\right)
$$

Here $g_k(\MOM)=\pi k+f_k(|\MOM|)$ and $f_k$ is a smooth function that increases monotonically from 0 to $\pi k$ on
an interval that will be specified later. Outside this interval, $f_k$ will
be identically equal to 0 or $\pi k$. Though the details do not matter for the Dehn twist
itself, our computations will turn out to put some constraints on $f_k$. 

For small $|\MOM|$, the map $\Dehn_k$ equals $(-1)^k\id$, while for large
$|\MOM|$ it equals the identity map.

\begin{defi}
  The map $\tau_k$ ($k\in\N$) is called a \textbf{$k$-fold right-handed Dehn twist}.
  The map $\tau_{-k}$ is called a \textbf{$k$-fold left-handed Dehn twist}.
\end{defi}

We will now construct a mapping torus of $T^*\S^{n-1}$ using these Dehn
twists following the construction of Giroux and Mohsen
\cite{Giroux_talk}. The canonical $1$-form $\lcan=\MOM \cdot d \POS$ on $T^*\S^{n-1}$ transforms like
$$ \Dehn_k^* \lcan = \lcan + |\MOM|\,d\big(f_k(|\MOM|)\big). $$
Note that the difference $ \lcan - \Dehn_k^* \lcan $ is exact. This implies
in particular that the Dehn twists are symplectomorphisms of $(T^*\S^{n-1},
d\lcan)$.
As a primitive of this difference $ \lcan - \Dehn_k^* \lcan $ we take 
$$
h_k(|\MOM|):=1-\int_0^{|\MOM|} s
f_k'(s)ds.
$$ 
Note that $h_k$ can assumed to be positive by choosing a
suitable interval where $f_k$ increases. To be more explicit, choose a
smooth function $f$ that is identically 0 on the interval $[0,1]$, on the interval
$[1,2]$ it increases monotonically from 0 to 1 and $f$ is identically 1 on the interval
$[2,\infty)$. Furthermore, we may assume that the derivative $f'$ is bounded
by 2. Then we can take $f_k(x):=k\pi f(c_k x)$ with $c_k>3 k\pi$. We have
$$
\int_0^{|p|} s f_k'(s) ds \leq \int_0^\infty k\pi c_k s f'(c_k s) ds \leq
k \pi \int_0 ^\infty y f'(y) dy/c_k \leq \frac{k \pi}{c_k} \int_1^2 y 2 dy
=\frac{3k \pi}{c_k},
$$
where we have substituted $y=c_k s$ and used that $f'(y)=0$ outside the
interval $[1,2]$ and that $f'$ is bounded by 2. Our choice of $c_k$ insures
that this integral is indeed smaller than 1, so $h_k$ is positive.
Consider the map
\begin{eqnarray*}
\phi_k:~ \R \times T^*\S^{n-1} & \longrightarrow & \R \times T^*\S^{n-1},\\
(t;\POS,\MOM) & \longmapsto & (t+h_k(|\MOM |);\tau_k(\POS,\MOM)).
\end{eqnarray*}
This map preserves the contact form $dt+ \lcan$ on $\R \times T^*\S^{n-1}$,
so we obtain an induced contact structure on  $\R \times T^*\S^{n-1} / \phi_k$.
\par
To make computations more convenient, we construct an additional
intermediate mapping torus. Let $\R \times T^*\S^{n-1}/\sim_k$ be the mapping
torus obtained by identifying $(t;\POS,\MOM) \sim_k (t+1;\tau_k(\POS,\MOM))$. 
We can define a diffeomorphism 
$$\R \times T^*\S^{n-1} /\sim_k \to \R \times T^*\S^{n-1} /\phi_k$$
by sending $(t;\POS,\MOM)$ to $(h_k(|\MOM |) t;\POS,\MOM)$. The pull-back
$\beta_k$ of the
described contact form under this diffeomorphism is given by
$$\beta_k=h_k(|\MOM |)dt-t|\MOM|\,d\big(f_k(|\MOM|)\big)+\lcan$$.

\section{Open books for the Brieskorn manifolds $\Brieskorn{2n-1}{k}$}\label{brieskorn}

The \textbf{Brieskorn manifolds} $\Brieskorn{2n-1}{k}\subset\C^{n+1}$ (with $k\in\N_0$) are defined as the intersection of the sphere $\S^{2n+1}$ with the zero set of the polynomial $f(z_0,z_1,\ldots,z_n)= z_0^k+z_1^2+\cdots+z_n^2$.
To make computations easier, assume that the radius of the $(2n+1)$-sphere
is $\sqrt 2$.

The orthogonal group $\SO(n)$ acts linearly on $\C^{n+1}$ by leaving the
first coordinate of $(z_0,z_1,\ldots,z_n)$ fixed and multiplying the last
$n$ coordinates with $\SO(n)$ in its standard matrix representation, i.e. $A\cdot(z_0,z_1,\ldots,z_n) := (z_0,A\cdot(z_1,\ldots,z_n))$.
This action restricts to $\Brieskorn{2n-1}{k}$, because the polynomial $f$ can be written as $z_0^k + \|\mathbf{x}\|^2 - \|\mathbf{y}\|^2+2i\pairing{\mathbf{x}}{\mathbf{y}}$ with $\mathbf{x}=(x_1,\ldots,x_n)$ and $\mathbf{y}=(y_1,\ldots,y_n)$.
  
Finally, the $\SO(n)$-invariant $1$-form
$$ \alpha_k := k\cdot(x_0\,dy_0-y_0\,dx_0) + 2\sum_{j=1}^n \left(x_j\,dy_j-y_j\,dx_j\right) $$
is of contact type on $\Brieskorn{2n-1}{k}$ for all $k\in\N$ as was shown
by Lutz and Meckert \cite{LutzMeckert}.

It is well-known that all $\Brieskorn{2n-1}{k}$ are $(n-2)$-connected and some
of these Brieskorn manifolds are spheres \cite{Brieskorn}, \cite{Hirzebruch}.
Ustilovsky \cite{Ustilovsky} showed that among them there are
diffeomorphic but non-contactomorphic manifolds. Namely if $2n-1=1
~\mathrm{mod}~4$, then all $\Brieskorn{2n-1}{k}$ with $k=\pm
1~\mathrm{mod}~8$ are standard spheres with inequivalent contact structures.
\par 
In the remainder of this paper will we show that the contact structures on
Brieskorn manifolds $\Brieskorn{2n-1}{k}$ are supported by an open book
whose monodromy is given by a $k$-fold Dehn twist. We define the binding
$B$ of the open book by the set in $\Brieskorn{2n-1}{k}$ with
$z_0=0$. We have the fibration $\theta:\,(\Brieskorn{2n-1}{k}-B)\to\S^1$,
given by $(z_0,z_1,\dots,z_n)\mapsto z_0/|z_0|$.

\subsection{The binding}\label{binding}
The only stabilizers of the $\SO(n)$-action on the Brieskorn manifold that occur are $\SO(n-1)$ and $\SO(n-2)$.
The projection onto the orbit space is given by
\begin{eqnarray*}
\Brieskorn{2n-1}{k} & \longrightarrow & \D^2 \\
(z_0,z_1,\ldots,z_n)& \longmapsto  & z_0.
\end{eqnarray*}
Points $(z_0,\dots,z_n)$ lying over the interior of the disk (i.e. $|z_0|\ne 1$) have principal stabilizer, points over $\partial\D^2$ lie on singular orbits.
The orbit $B=\orb{(0,1,i,0,\dots,0)}\cong \SO(n)/\SO(n-2)$ is the binding of the open book.
It is naturally contactomorphic to ${W_2^{2n-3}}$. In fact,
$W_2^{2n-3}=\SO(n)/\SO(n-2)$ is diffeomorphic to the unit sphere bundle
$\S(T^*\S^n)$. This shows that part (1) of Definition~\ref{compatiblebook} is satisfied.

The symplectic normal bundle of the binding is trivial, because for $k \neq 1$ we have a symplectic basis
$$\frac{1}{\sqrt {2k}}(1,0,\ldots,0)\text{,~~~}\frac{1}{\sqrt
  {2k}}(i,0,\ldots,0),$$
and for $k=1$ we have the basis
$$\sqrt{\frac{2}{5}}(1,-\frac{\bar z_1}{4},\ldots,-\frac{\bar
  z_n}{4})\text{,~~~}\sqrt{\frac{2}{5}}(i,-\frac{i\bar
  z_1}{4},\ldots,-\frac{i\bar z_n}{4}).$$ 
The neighborhood theorem for contact submanifolds \cite{Geiges_Preprint}
then shows that there is a neighborhood of the binding that is
contactomorphic to $(B\times \D^2,\alpha_k|_B+r^2 d\theta)$, where
$(r,\theta)$ are polar coordinates on the disk.

\subsection{The pages}

In this section, we will prove that $\Brieskorn{2n-1}{k}-B$ is
contactomorphic to $\R\times T^*\S^{n-1}/\sim_k$, the mapping torus of a
$k$-fold Dehn twist. 

The $\R$-action on $\Brieskorn{2n-1}{k}-B$, given by 
$$
e^{it}(z_0,z_1,\dots,z_n)=(e^{it}z_0,e^{\frac{ki}{2}t}z_1,\ldots,e^{\frac{ki}{2}t}z_n).
$$
induces a diffeomorphism between the pages $\theta^{-1}(1)$ and $\theta^{-1}(e^{it})$.

Let us define an auxiliary mapping torus to make computations more
convenient. Define 
$$M_k:=\R\times T^*\S^{n-1}/\sigma_k,$$ 
where 
$$\sigma_k(t,\POS,\MOM)=(t+1,(-1)^k\POS,(-1)^k\MOM).$$ We will now give an explicit map to show that
$P$ is diffeomorphic to $T^*_{|\MOM|<1}\S^{n-1}$. Here $T^*_{|\MOM|<1}\S^{n-1}$ denotes the
open unit disk bundle associated with the cotangent bundle of $\S^{n-1}$.
A point $(\POS,\MOM)\in T^*\S^{n-1}\subset \R^n\times \R^n$ with $|\POS|=1$,
$|\MOM|\le 1$, and $\POS\perp\MOM$ is mapped to
$$ (\POS,\MOM) \mapsto \Big(1-|\MOM|^2,F(|\MOM|)\MOM + i G(|\MOM|)\POS \Big) $$
with $F(r)=\sqrt{\frac{2-(1-r^2)^2-(1-r^2)^k}{2r^2}}$ and $G(r)=\sqrt{\frac{2-(1-r^2)^2+(1-r^2)^k}{2}}$.

Together with the $\R$-action this gives a map
\begin{eqnarray*} 
\Phi_k:\,\R\times T^*_{|\MOM|<1}\S^{n-1} & \longrightarrow & \Brieskorn{2n-1}{k}\\
(t,\POS,\MOM) & \longmapsto & \Big(e^{2\pi it}(1-|\MOM|^2),e^{\pi kit}(F(|\MOM|)\MOM + i
G(|\MOM|)\POS) \Big). 
\end{eqnarray*}
This descends to a diffeomorphism of the subset of $M_k$ with $|\MOM|<1$ to
$\Brieskorn{2n-1}{k}-B$. For $k$ even, one obtains
$\Phi_k(t+1,\POS,\MOM)=\Phi_k(t,\POS,\MOM)$, so that
$\Brieskorn{2n-1}{k}-B\cong\S^1\times T^*_{|\MOM|<1}\S^{n-1}$, and for $k$ odd, one
obtains $\Phi_k(t+1,\POS,\MOM)=\Phi_k(t,-\POS,-\MOM)$, so that
$\Brieskorn{2n-1}{k}-B$ is a non-trivial $T^*_{|\MOM|<1}\S^{n-1}$-bundle over $\S^1$.

The pull-back of the contact form $\alpha_k$ to $M_k$ under $\Phi_k$ gives
$$ \Phi_k^*\alpha_k =2\pi k\big((1-|\MOM|^2)^2 + |\MOM|^2F^2 + G^2\big)\,dt + 4FG\,\lcan = 4\pi k\,dt + 4FG\,\lcan. $$

Next, we construct a diffeomorphism $\Psi_k$ from $M_k$ to the mapping
torus $\R\times T^*\S^{n-1}/\sim_k$ by defining
$$  
   \Psi_{k}(t;\POS,\MOM)= \Big[t;\POS\cdot\cos\big(tf_k(|\MOM|)\big) + \frac{\MOM}{|\MOM|}\cdot\sin\big(tf_k(|\MOM|)\big), 
   \MOM\cdot\cos\big(tf_k(|\MOM|)\big) - |\MOM|\POS\cdot\sin\big(tf_k(|\MOM|)\big)\Big].
$$
The map is well-defined, because $\Psi_k \circ \sigma_k(t;\POS,\MOM)$ is
identified with $\Psi_k(t;\POS,\MOM) $ in the mapping torus $\R\times T^*\S^{n-1}/\sim_k$. 
In order to show that $(\Brieskorn{2n-1}{k}-B,\alpha_k)$ and $(\R \times
T^*\S^{n-1} /\sim_k,\beta_k)$ are contactomorphic, we will show that the
pull-back of $\alpha_k$ under $\Phi_k$ is contactomorphic to the pull-back of $\beta_k$ under $\Psi_k$. 

We now compute the pull-back of $\beta_k$ under $\Psi_k$, noting that the
norm of $p$ is invariant under $\Psi_k$ (we do not write the dependence of
$h_k$ and $f_k$ on $|\MOM|$):
\begin{multline*}
\Psi_k^* \beta_k  = h_k dt-t|\MOM|d f_k+\left( \MOM \cos(tf_k)-|\MOM|\POS\sin(tf_k)
\right) \cdot \Big( d\POS \cos(tf_k)- \\
-\POS \sin(tf_k)(f_k dt+ t d f_k)+  (\frac{d\MOM}{|\MOM|}-\frac{\MOM d|\MOM|}{|\MOM|^2})\sin(tf_k)+\frac{\MOM}{|\MOM|}\cos(tf_k)(f_k
dt+t df_k) \Big).
\end{multline*}
Since we have $\MOM \cdot \POS=0$ and $|\POS|^2=1$, it follows that $\MOM d
\POS=-\POS d \MOM$ (recall that $\MOM d\POS=\lcan$) and $\POS d \POS=0$.
We now use the standard trigonometric equalities and the fact that
$h_k(y)=1-yf_k(y)+\int_0^yf_k(s)ds$ to find
$$
\Psi_k^* \beta_k=\left( 1+\int_0^{|\MOM|}f_k(s)ds \right)\,dt + \lcan .
$$
Note that $\Phi_k^*\alpha_k$ has a very similar form.
We make the following ansatz for a contactomorphism of $(M_k|_{|\MOM|<1},\Phi_k^*\alpha_k)$ to $(M_k,\Psi_k^*\beta_k)$:
\begin{eqnarray*}
 S_k: (t,\POS,\MOM) & \mapsto & (t,\POS,\frac{g(|\MOM|)}{|\MOM|} \MOM).
\end{eqnarray*}
With this ansatz we find what $\MOM$ should map to in order to be a
contactomorphism. Note that we just rescale $\MOM$. The pull-back under this map of $\Psi_k^* \beta_k$ is given by
$$ 
\left( 1+\int_0^{g(|\MOM|)}f_k(s)ds \right) dt+\frac{g(|\MOM|)}{|\MOM|}\lcan.
$$
Since we want this to be a multiple of $\Phi_k^*\alpha_k$ (we actually even
seek equality), we need to solve the following equation:
$$
\frac{g(|\MOM|)}{1+\int_0^{g(|\MOM|)}f_k(s)ds}=\frac{|\MOM|FG}{k\pi}.
$$
Define an auxiliary function 
$$
h(y):=\frac{y}{1+\int_0^y f_k(s)ds}.
$$ 
The above equation becomes
$$
h\left( g(|\MOM|) \right)=\frac{|\MOM|FG}{k\pi}.
$$
We will solve for $g(|\MOM|)$ by inverting $h$. This can be done by the
following considerations. The derivative of $h$ is given by 
$$
h'(y)=\frac{1-\int_0^y s f_k'(s)ds}{\left( 1+\int_0^y f_k(s)ds \right)^2}
=\frac{h_k(y)}{\left( 1+\int_0^y f_k(s)ds \right)^2}
$$ and is positive by our choice of $h_k$
in Section~\ref{dehn twist}. Since this shows that $h$ is strictly increasing, we also observe
that the function $h$ maps $[0,\infty)$ to $[0,\frac{1}{k \pi})$. This can be
seen by noting that $f_k(s)=k \pi$ for $s$ sufficiently large, again due to our
choice of $h_k$. It also means
that $h$ can be inverted when restricted to a suitable range. One easily
checks that the right-hand side of the above equation,
$\frac{|\MOM|FG}{k\pi}$, has positive derivative and is therefore strictly
increasing on the interval $[0,1)$. Moreover it has the same
range as $h$, namely $[0,\frac{1}{k \pi})$. Therefore we can find a smooth
solution to $g(|\MOM|)$ by applying the inverse of $h$ to $\frac{|\MOM|FG}{k\pi}$.

This shows that the open book $(B,\theta)$ on $\Brieskorn{2n-1}{k}$ has
page $T^*\S^{n-1}$ with monodromy given by a $k$-fold Dehn twist. The
contactomorphism that achieves this is
$$C_k:=\Phi_k \circ S_k^{-1} \circ \Psi_k^{-1}: (\R\times T^*\S^{n-1}/\sim_k,\beta_k) \to
(\Brieskorn{2n-1}{k}-B,\alpha_k).
$$
Note that this contactomorphism also respects the projection to $\S^1$,
because the $\S^1$-coordinate is invariant under $C_k$.

\subsection{The contact structure on $\Brieskorn{2n-1}{k}$ is supported by the open book.}
Part (1) of the Definition~\ref{compatiblebook} was already checked in Section~\ref{binding}. With our
identification of the pages with $T^*\S^{n-1}$ and the fact that $C_k$ respects the fibration, the other two parts follow
easily. We have that $d\beta_k$ restricts to $d\lcan$ on
any page, so $(T^*\S^{n-1},d\beta_k)$ is symplectic and therefore any page of
$\theta$ on $\Brieskorn{2n-1}{k}-B$ will also be
symplectic with form $d\alpha_k$. This shows part (2).

Instead of showing property (3) directly, take a page $P$ and find a copy
of the binding in $P$. We do the latter as follows. Take $w_0\in \C-\{0\}$ such
that $P=\theta^{-1}(w_0/|w_0|)$. Then define $B_{w_0}:=\{ (z_0,z_1,..,z_n)\in
\Brieskorn{2n-1}{k} ~|~ z_0=w_0\}$. Note that $B_{w_0}$ is contactomorphic
to $B$ if $|w_0|$ is small enough and has therefore the same orientation induced by its contact form.

We can now regard $B_{w_0}$ as the boundary of the compact page $P_{w_0}$ by
cutting off a part of $P$. The computations are a little easier if we
consider the inverse images of $B_{w_0}$ and $P_{w_0}$ under the map
$C_k$.

The inverse image of $B_{w_0}$ under $C_k$ is the set 
$$\left\{ \frac{w_0}{|w_0|} \right\}\times \S_cT^*\S^{n-1}\subset \R\times T^*\S^{n-1}/\sim_k,$$
where $\S_cT^*\S^{n-1}$ denotes the associated circle bundle with radius $c$
and the constant $c$ is determined by $|w_0|$. Similarly, the trimmed page $P_{w_0}$ is mapped to
$$\left\{  \frac{w_0}{|w_0|} \right\} \times T_{|\MOM|\leq c}^*\S^{n-1}.$$ This means that it suffices to check that
the orientation of $(\S_cT^*\S^{n-1},\lcan)$ induced by the contact form
coincides with its orientation as the
boundary of $(T_{|\MOM|\leq c}^*\S^{n-1},d\lcan)$. This is done by observing
that $\MOM \frac{\partial}{\partial \MOM}$ is a Liouville vector field for $d\lcan$. 

\bibliography{openbooks}

\providecommand{\bysame}{\leavevmode\hbox to3em{\hrulefill}\thinspace}
\providecommand{\MR}{\relax\ifhmode\unskip\space\fi MR }
\providecommand{\MRhref}[2]{%
  \href{http://www.ams.org/mathscinet-getitem?mr=#1}{#2}
}
\providecommand{\href}[2]{#2}
\begin{thebibliography}{HM68}

\bibitem[Bri66]{Brieskorn}
E.~Brieskorn, \emph{Beispiele zur {D}ifferentialtopologie von
  {S}ingularit\"aten}, Invent. Math. \textbf{2} (1966), 1--14. \MR{34 \#6788}

\bibitem[Gei]{Geiges_Preprint}
H.~Geiges, \emph{Contact geometry}, to appear in the Handbook of Differential
  Geometry, vol. 2, Elsevier, arXiv:math.SG/0307242.

\bibitem[Gir02]{Giroux}
E.~Giroux, \emph{G\'eom\'etrie de contact: de la dimension trois vers les
  dimensions sup\'erieures}, Proceedings of the International Congress of
  Mathematicians, Vol. II (Beijing), Higher Ed. Press, 2002, pp.~405--414.
  \MR{2004c:53144}

\bibitem[GM]{Giroux_talk}
E.~Giroux and J-P. Mohsen, \emph{Contact structures and symplectic fibrations
  over the circle, lecture notes}.

\bibitem[HM68]{Hirzebruch}
F.~Hirzebruch and K.~H. Mayer, \emph{{${\rm O}(n)$}-{M}annigfaltigkeiten,
  exotische {S}ph\"aren und {S}ingularit\"aten}, Lecture Notes in Mathematics,
  No. 57, Springer-Verlag, Berlin, 1968. \MR{37 \#4825}

\bibitem[LM76]{LutzMeckert}
R.~Lutz and C.~Meckert, \emph{Structures de contact sur certaines sph\`eres
  exotiques}, C. R. Acad. Sci. Paris S\'er. A-B \textbf{282} (1976), no.~11,
  Aii, A591--A593. \MR{53 \#1471}

\bibitem[Sei]{Seidel_Preprint}
P.~Seidel, \emph{{Symplectic automorphisms of $T^*S^2$}},
  arXiv:math.DG/9803084.

\bibitem[Ust99]{Ustilovsky}
I.~Ustilovsky, \emph{Infinitely many contact structures on {$S\sp {4m+1}$}},
  Internat. Math. Res. Notices (1999), no.~14, 781--791. \MR{2000f:57028}

\end{thebibliography}

\end{document}